\documentclass[11pt]{article}
\usepackage{geometry}                
\usepackage{graphicx}
\usepackage{amssymb}
\usepackage{epstopdf}
\newtheorem{theorem}{Theorem}
\newtheorem{proposition}{Proposition}
\newtheorem{example}{Example}
\newtheorem{definition}{Definition}

\newtheorem{lemma}{Lemma}
\newtheorem{corollary}{Corollary}

\def\R{\Bbb R}
\def\C{\Bbb C}

\def\H{\Bbb H}
\def\Z{\Bbb Z}
\def\spp{\vspace{5pt}  \noindent}

\newcommand{\eqref}[1]{{\rm (\ref{#1})}}

\input mssymb.tex

\title{Conditional Extremals}
\author{Lyle Noakes\\
School of Mathematics and Statistics,\\
The University of Western Australia,\\
Nedlands WA 6009,\\
Australia.\\
{\tt lyle@maths.uwa.edu.au}}

\begin{document}
\maketitle
\begin{center}{\Large Abstract} 
\end{center}
Imagine that measurements are made at times $t_0$ and $t_1$ of the trajectory of a physical system whose governing laws are given approximately by a class ${\cal A}$ of so-called {\em prior vector fields}. Because the physical laws are not known precisely, it might be that the measurements are not realised by the integral curve of any prior field. We want to estimate the behaviour of the physical system between times $t_0$ and $t_1$.  
  
 \spp
An integral curve of an arbitrary vector field $X$ is said to be {\em feasible} when it interpolates the measurements. When $X$ is critical for $L^2$ distance to ${\cal A}$, the feasible curve is called a {\em conditional extremum}. When the distance to ${\cal A}$ is actually minimal, the conditional extremum is a best estimate for the intermediate behaviour of the system. 

\spp
The present paper does some of basic groundwork for investigating mathematical properties of conditional extrema,  focusing on cases where ${\cal A}$ comprises a single prior field $A$. When $A={\bf 0}$ a conditional extremal is a geodesic arc, but this special case is not very representative. 
In general, $A$ enters into the Euler-Lagrange equation for conditional extrema, and more can be said when $A$ is conservative or has special symmetry. 

\spp
We characterise conservative priors on simply-connected Riemannian manifolds in terms of their conditional extrema:  when time is reversed, a constant is added to the $L^2$ distance. For some symmetric priors on space forms we obtain conditional extrema in terms of the Weierstrass elliptic function. 
For left-invariant priors on bi-invariant Lie groups, conditional extrema are shown to be right translations of pointwise-products of $1$-parameter subgroups.  
\section{Introduction} This paper focuses on the following question. 

\spp
{\em Let $A$ be a given $C^\infty$ vector field on a connected $C^\infty$ $m$-manifold $M$. 
A curve $x$ in $M$ is thought to be an integral curve of a unknown vector field that is near $A$. Precise observations $x(t_k)$ are made of $x$ 
at finitely many times $t_k$ where  
$t_0\leq t_1\leq \ldots \leq t_n$, but no further information about $x$ is known. How then to estimate $x(t)$ for $t\not= t_k$?}  

\spp
Consider first the case where $M$ is Euclidean $m$-space $E^m$. If $n=1$ and $t\in (t_0,t_1)$ we could use the weighted average 
$$x(t)~\approx ~\frac{(t_1-t)x_0+(t-t_0)x_1}{t_1-t_0}$$
or, more generally, for $n>1$ the natural cubic spline \cite{deboor}.  But these interpolants take no account of the information that 
$x^{(1)}\approx A(x)$, where $x^{(1)}$ denotes the derivative of $x$. To do so requires  
knowledge about the nature of the supposed approximation, which in practice depends on context. 
We are going to assume $x^{(1)}$ is close to $A(x)$ in the $L^2$ sense, namely 
$x^{(1)}(t)=A(x(t))+g(t)$ where $g:[t_0,t_n]\rightarrow E^m$ is $L^2$-small. Possibly $g(t)=G(x(t))$ where $G$ is another vector field.  

\spp
In practice there might be a class ${\cal A}$ of vector fields rather than a single prior field $A$, with the information that $x^{(1)}\approx A(x)$ for some $A\in {\cal A}$. When ${\cal A}$ is large the quality of the information about the field is small. For ${\cal A}$ the space of all $C^\infty$ vector fields there is no information except the observations, and we are reduced to classical methods of interpolation. Between these extremes, where ${\cal A}$ is a singleton or very large, 
${\cal A}$ might instead be parameterised by a finite-dimensional $\alpha \in \R ^p$. 
The literature on model building and parameter selection for dynamical systems (see  \cite{nelles},  \cite{chu}) includes cases where observations are contaminated by noise, as in \cite{toni}. In the simpler setting of our question, there is no noise, and the requirement is for interpolation rather than approximation. An important case for the present paper is where ${\cal A}$ is a singleton, and parameter selection is not an issue. 

\spp 
Our task is to study nonparametric interpolants $x$ minimising the $L^2$ difference between $x^{(1)}$ and $A(x)$ subject only to $x(t_k)=x_k$ for $0\leq k\leq n$ - a problem in the calculus of variations.  The present work also finds new results and explicit solutions for systems with symmetry in the setting of Riemannian geometry, where $E^m$ is replaced by an $m$-dimensional Riemannian manifold $M$. 
\section{Description of Results}
In \S \ref{optsec} conditional minima and extrema $x$ are defined relative to a class ${\cal A}$ of vector fields. Working first within the class of so-called almost-$C^2$ curves, it follows from Theorem \ref{thm1} in  \S \ref{elsec}, that if $A$ is $C^k$ then, except at $t_1,t_2,\ldots ,t_{n-1}$, a conditional extremum $x$ is automatically $C^{k+1}$. Examples \ref{ex0}, 
\ref{geoex}, \ref{eucex1} point out some simple facts. For instance, when ${\cal A}=\{ {\bf 0}\}$ a conditional extremum is the same as a geodesic: conditional extrema relative to prior fields generalise geodesics. 

\spp
Theorem \ref{thm1} also gives the Euler-Lagrange equation (\ref{eleq}) for a conditional extremum $x$ with respect to a prior field $A$. The left hand side of (\ref{eleq}) is the covariant acceleration of $x$. The right hand side has two nontrivial terms, namely the gradient of the squared norm of $A$, and a second term made by 
\begin{itemize}
\item replacing $A$ by the exterior $1$-form $A^T$ corresponding under the Riemannian metric
\item taking the exterior derivative of $A^T$
\item contracting with $x^{(1)}$
\item replacing the resulting $1$-form by the corresponding vector field defined along $x$. 
\end{itemize}
A conditional extremum is a $C^0$ track-sum of solutions of (\ref{eleq}). Consequently, when ${\cal A}$ is a singleton  it suffices to consider the case $n=1$. For $M=E^m$, Example \ref{ex1}  solves equation (\ref{eleq}) when $A$ is an affine vector field (Example \ref{eucex1} finds $x$ directly when $A$ is constant). Theorem \ref{cor0} shows, for an arbitrary Riemannian manifold $M$, that $\Vert x^{(1)}(t)\Vert ^2-\Vert A(x(t))\Vert ^2$ is conserved along solutions $x$ of (\ref{eleq}), and that 
if $A$ is bounded and $M$ is complete then $x$ extends to a solution of (\ref{eleq}) defined on the whole of $\R $. 
This generalises the well-known results for geodesics, that $\Vert x^{(1)}\Vert ^2$ is conserved, and that complete Riemannian manifolds are geodesically complete \cite{milnor} Part II. 

\spp
For a geodesic $x:[t_0,t_1]\rightarrow M$ its {\em reverse} $\bar x:[-t_1,-t_0]\rightarrow M$, given by $\bar x(u)=x(-u)$ is also a geodesic. From the asymmetric role of the prior field $A$ in the definition, it seems improbable that the reverse of a conditional extremum would be a conditional extremum for the time-reversed data. The improbable happens when the $1$-form $A^T$ is closed, as proved in Corollary \ref{prop1} of Theorem \ref{thm1}. The reverse of a conditional optimum need not be optimal, even when $A^T$ is closed, as seen in Example \ref{counterprop1}. 

\spp
A stronger condition is for $A$ to be {\em reflexive}, namely $J(x,A)$ differs from $\bar J(\bar x,A)$ by a constant depending only on $A$ and independent of $x$. Theorem \ref{thm2} of \S \ref{closedsec} shows $A$ is reflexive if $A^T$ is exact and, conversely, if $A$ is reflexive then $A^T$ is closed. So if $M$ is simply connected then reflexive is equivalent to conservative, as illustrated in Example \ref{refex} by numerical calculation. These kinds of calculations are performed by replacing curves $x$ by lists of points in $\R ^p$ where $p$ is $2$, $3$ or $4$, then numerically optimising with Mathematica's {\em FindMinimum}.  In applications one might rely on refinements of such methods, but the present paper gives theoretical results, including closed form solutions for conditional optima, when $M$ and $A$ exhibit symmetry. 

\spp
In \S \ref{Liesec}, $M$ is a semisimple Lie group $G$ with bi-invariant Riemannian metric. Whereas geodesics in $G$ are right-translations of $1$-parameter subgroups,  Theorem \ref{thm2g} says solutions of (\ref{eleq}) for a left-invariant prior field $A$ are right-translations of {\em pointwise products} of a {\em pair} of $1$-parameter subgroups. One subgroup is generated by the value of $A$ at the identity of $G$. The other is obtained from the first after comparison with $x_0$ and $x_1$. Corollary \ref{cor11} extends this to $n>1$. So the Euler-Lagrange equation (\ref{eleq}) is solved in closed form, as illustrated for the group $S^3\cong SU(2)$ of unit quaternions  by Example \ref{S3ex}. All this supposes that the left-invariant field $A$ is given, as it would be if ${\cal A}$ was a singleton. For a finite set of left-invariant fields, $J$ is minimised for each $A$, and the smallest is chosen. 

\spp
For ${\cal A}$ an infinite set of left-invariant vector fields on $G$, Corollary \ref{cor11} of Theorem \ref{thm2g} reduces the search for an optimal pair to a finite-dimensional optimisation problem, whose numerical solution is straightforward. We also prove some theoretical results for the case when ${\cal A}$ is generated by a submanifold of ${\cal G}$. Theorem \ref{sumBthm} gives a nonlinear equation for $A$ in terms of the exponential map of $G$.  When ${\cal A}$ is the set of all left-invariant fields and the $t_k$ are equally spaced, the equation takes a simpler form given in Corollary \ref{corsum}, and verified in Example \ref{S3ex2}. Specialising to $n=2$ (three observations), Corollary \ref{3pts} gives a simple solution for $(x,A)$ in terms of the exponential.

\spp
From \S \ref{sphsec} onwards, $M$ is either the unit sphere $S^2$ in $E^3$ or the unit two-sheeted hyperboloid $H^2$ in Lorentz $3$-space, the two cases being treated simultaneously. The prior fields $A$, parameterised by a pair of functions $\bar \beta ,\bar \gamma :\R \rightarrow \R$, are invariant with respect to rotations in $\R ^2\times \{ 0\} $. Under these conditions, in \S \ref{sphsec} equation (\ref{eleq}) is rewritten  as three coupled nonlinear $2$nd order scalar ODEs (\ref{x1eq}), (\ref{x2eq}), (\ref{x3eq}) for the scalar-valued coordinates $x_1$, $x_2$, $x_3$ of a conditional extremum $x$.  As well as equation (\ref{cons}) which reappears as (\ref{conssph}), rotational symmetry gives another conserved quantity (\ref{consrot}). 
Then $x$ is found by quadrature in terms of $x_3:\R \rightarrow \R$, and a first order ODE (\ref{x3eqss}) for $x_3$ is given. Then $\bar \beta$ and $\bar \gamma $ are taken as constant, and the solutions of (\ref{x1eq}), (\ref{x2eq}), (\ref{x3eq}) for which $x_3$ is constant are found. Further results for $x_3$ nonconstant depend on whether $A$ has a longitudinal component, namely whether $\bar \gamma \not= 0$. 

\spp
For $\bar \gamma =0$, Theorem \ref{horthm} of \S \ref{sphsec} solves (\ref{x1eq}), (\ref{x2eq}), (\ref{x3eq}) for $x$ 
in terms of $\sin $, $\cos $, $\sinh $, $\arctan $, $\cosh $ and ${\rm arctanh} $. For both $S^2$ and $H^2$ there are two different kinds of solution, depending on whether $(0,0,1)$ lies in the image of $x$. Examples \ref{hor2ex}, \ref{poin1ex} describe solutions where $M$ is $S^2$ and $H^2$ respectively, with $H^2$ replaced by the Poincar\'e unit disc for ease of illustration. As shown in Figure \ref{horconfig}, there exist non-optimal conditional extrema. 

\spp
In \S \ref{gamnon0sec}, $\bar \gamma \not= 0$ and Theorem \ref{wthm} gives $x_3$ in terms of the Weierstrass elliptic function $\wp$. In particular $x_3$ is periodic, and solutions $x$ of (\ref{x1eq}), (\ref{x2eq}), (\ref{x3eq}) are found by quadrature in terms of $\wp$. Even when $A$ is conservative, namely $\bar \beta =0$, these solutions can be geometrically interesting as shown for $M=S^2$ in Figure \ref{fig0}. When $A$ is not conservative $A$ is usually not reflexive (never on $S^2$ or $H^2$), and Example \ref{poin2ex} describes a conditional extremum $x:[0,1]\rightarrow H^2$ and its reverse $\bar x$ where $A$ is nonconservative.  This is illustrated in Figure \ref{poin2fig} by mapping into 
the Poincar\'e disc.   
\section{Optimality}\label{optsec}
Let $\langle ~,~\rangle $ be a $C^\infty$ Riemannian metric on $M$. For $a_0<a_1\in \R$, a continuous  curve $y:[a_0,a_1]\rightarrow M$ is said to be {\em almost-}$C^2$ when, for some $\alpha >0$, some $a_0=s_0<s_1<\ldots <s_p=a_1$ and all $j=1,2,\ldots p$, each restriction $y\vert [s_{j-1},s_j]$ extends to a $C^{2}$ curve in $M$ defined on $(s_{j-1}-\delta ,s_j+\delta )$. Then $s_0,s_1,\ldots s_{p}$ are {\em singular parameters} of $y$. 
%
%
%
%

\spp 
For some $n\geq 1$, let $x_0,x_1, \ldots x_n\in M$ and $t_0<t_1<\ldots <t_n\in \R$ be given. An almost-$C^{2}$ curve $x:[t_0,t_n]\rightarrow M$ is said to be {\em feasible} when $x(t_k)=x_k$ for all $k=1,2,\ldots ,n$. 
The set of all feasible curves is denoted by 
${\cal X}={\cal X}_{t_0,t_1,\ldots ,t_n;x_0,x_1,\ldots x_n}$. 
Given also a nonempty set ${\cal A}$ of $C^{1}$ vector fields $A$ on $M$, define $J:{\cal X}\times {\cal A}\rightarrow \R $ by 
$$J(x,A)~:=~\int _{t_0}^{t_n}\Vert x^{(1)}(t)-A(x(t))\Vert  ^2~dt$$ 
where $x^{(1)}$ is the derivative of $x$ with respect to $t$, and $\Vert ~\Vert $ denotes the Riemannian norm. 
A minimizer $(x,A)$ of 
$J=J_{t_0,t_1,\ldots ,t_n;x_0,x_1,\ldots x_n}$ is said to be {\em conditionally-optimal} or just {\em optimal}. When there is no doubt about $A$, for instance when ${\cal A}$ is a singleton, we say $x$ is conditionally optimal rather than $(x,A)$. 

\spp
Suppose there is an unknown curve $z:[t_0,t_n]\rightarrow M$ approximately satisfying 
$$z^{(1)}(t)~=~A(z(t))$$ 
for some vector field $A$ in a given parameterized set ${\cal A}$ of {\em prior fields}. 
To say $(x,A)$ is conditionally optimal means the velocity field of the feasible curve $x$ is as near as possible to $A\in {\cal A}$ while agreeing with observations made from $z$ at the parameter values $t_0,t_1,\ldots ,t_n$. 
\begin{example}\label{ex0} If some $x\in {\cal X}$ is an integral curve for some $A\in {\cal A}$ then $(x,A)$ is conditionally optimal. For $n=1$ and $x_0=x_1$, the constant curve is not necessarily optimal. $\Box$
\end{example} 
\begin{example}\label{geoex} Let ${\cal A}=\{ {\bf 0}\} $. Then $x\in {\cal X}$ is conditionally optimal when it is a track 
sum of $n$ minimal geodesic arcs from the $x_{k-1}$ to $x_k$. In particular, when $n=1$, $x$ is a minimal geodesic 
arc. $\Box$
\end{example}
\begin{example}\label{eucex1} Let $M$ be Euclidean $m$-space $E^m$, and let $A$ be a constant field. Then 
$$J(x,A)~=~\int _{t_0}^{t_n}\Vert x^{(1)}(t)\Vert ^2~dt -2\langle A,x_n-x_0\rangle +(t_n-t_0)\Vert A\Vert ^2.$$
So, for any nonempty set ${\cal A}$ of constant vector fields on $E^m$, a necessary condition for $(x,A)$ to be conditionally optimal is that 
$x$ be the piecewise-affine interpolant of the observed points $x_k$ at the $t_k$. As will be seen, it is rare for $x$ to be independent of ${\cal A}$ in this way. 

\spp
In the least informative situation, where ${\cal A}$ is the set of all constant fields there is a unique optimal pair $(x,A)$, with  
$A=(x_n-x_0)/(t_n-t_0)$. It might seem strange that the optimal $A$ takes no account of $x_1,x_2,\ldots ,x_{n-1}$, 
but the observations are noise-free. The uncertainty concerns only $x(t)$ for $t\not= t_k$. So $A$ should indeed be a weighted mean of intermediate estimates:
$$A~=~\sum_{k=1}^n(t_k-t_{k-1})\left( \frac{x_{k}-x_{k-1}}{t_k-t_{k-1}}\right) .$$ 
The outcome is less predictable in Example \ref{ex1} where ${\cal A}$ is the set of all affine vector fields on $E^m$, and in  \S \ref{lsr} where ${\cal A}$ is all left-invariant vector fields on a semisimple Lie group. $\Box$
\end{example}
\section{First Order Necessary Conditions} \label{elsec}
Let $\tilde x:[t_0,t_n]\times (-\epsilon ,\epsilon )\rightarrow M$ be continuous where $\epsilon >0$.  Set $x_h(t):=x^t(h):=\tilde x(t,h)$. Given $x\in {\cal X}$ with singular parameters $s_0,s_1,\ldots ,s_p$, we call $\tilde x$ a {\em variation} of $x$ when  
\begin{itemize}
\item $\tilde x(t,0)=x(t)$ for all $t\in t_0,t_n]$
\item $x_h\in {\cal X}$ for all $h\in (-\epsilon ,\epsilon)$ 
\item for some $\alpha >0$, each restriction $\tilde x\vert [s_{j-1},s_j]\times (-\epsilon ,\epsilon )$ extends to a $C^{2}$ map from $(s_{j-1}-\alpha ,s_j+\alpha ) \times (-\epsilon ,\epsilon)$ to $M$. 
\end{itemize}
So  each $x_h$ has the same singular parameters as $x$, and $x^t$ is $C^{2}$ for all $t\in [t_0,t_n]$. Set $W(t)=(x^t)'(0)$ 
where $~'~$ denotes differentiation with respect to $h$ and, for $s\in (t_0,t_n)$,  
$\displaystyle{\Delta x^{(1)}(s):=\lim_{u\rightarrow s^+}x^{(1)}(u)-\lim_{u\rightarrow s^-}x^{(1)}(u)}$.
Then $\Delta x^{(1)}(s)={\bf 0}$ for $s$ nonsingular.
A vector field $X$ on $M$ corresponds  to a differential $1$-form $X^T$, given by $X^T(Y):=\langle X,Y\rangle $, where $Y$ is any vector field.  Of course $A$ is conservative precisely when $A^T$ is exact. Since $A\in {\cal A}$ $C^1$ so is $A^T$, and another differential $1$-form $\theta _{A,X}$ on $M$ 
is given by 
$$Y~\mapsto ~\theta _{A,X}(Y)~=~-\theta _{A,Y}(X)~:=~(dA^T)(X,Y)~=~$$
$$X(\langle A,Y\rangle )-Y(\langle A,X\rangle )-\langle A,[X,Y]\rangle ~=~\langle \nabla _XA,Y\rangle -\langle \nabla _YA,X\rangle $$
because the Levi-Civita covariant derivative $\nabla $ is torsion-free.   The
 differential $1$-form $\theta _{A,X}$ corresponds to another vector field 
$\theta _{A,X}^{-T}$. From the definition of $\theta _{A,X}$, $\theta _{A,X}(Y)=-\theta _{A,Y}(X)$ and, for any $C^1$ function $\alpha :M\rightarrow \R $,
\begin{equation}
\label{prodth}\theta _{\alpha A,X}~=~\alpha \theta _{A,X}+X(\alpha )A^T-\langle A,X\rangle d\alpha .
\end{equation}
\begin{lemma}\label{lem1} For any variation $\tilde x$ of $x\in {\cal X}$, 
$$\frac{1}{2}\frac{\partial }{\partial h}J(x_h,A)\vert _{h=0}~=~
\sum_{j=1}^{p-1}\langle W(s_j),\Delta x^{(1)}(s_j)\rangle -\int _{t_0}^{t_n}\langle W,\nabla _t x^{(1)}(t)-\frac{1}{2}{\rm grad }\Vert A\Vert ^2-\theta _{A,x^{(1)}}^{-T}\rangle 
~dt.$$ 
\end{lemma}

\spp
{\bf Proof:} Because $\nabla $ is symmetric and compatible with the Riemannian metric, 
$$\frac{1}{2}\frac{\partial }{\partial h}J(x_h,A)~=~
\int _{t_0}^{t_n}\langle \nabla _t x_h'-\nabla _hA,x_h^{(1)}(t)-A(x_h(t))\rangle ~dt$$
where $~'~$ denotes differentiation with respect to $h$. On integration by parts this becomes 
$$\sum_{j=1}^{p-1}\langle x_h'(s_j),\Delta x_h^{(1)}(s_j)\rangle -\int _{t_0}^{t_n}\langle x_h',\nabla _t (x_h^{(1)}(t)-A(x_h(t)))\rangle 
+\langle \nabla _hA,x_h^{(1)}(t)-A(x_h(t))\rangle 
~dt.$$
Setting ${h=0}$,   
$\displaystyle{\frac{1}{2}\frac{\partial }{\partial h}J(x_h,A)\vert _{h=0}~=}$
$$\sum_{j=1}^{p-1}\langle W(s_j),\Delta x^{(1)}(s_j)\rangle -\int _{t_0}^{t_n}\langle W(t),\nabla _t (x^{(1)}(t)-A(x(t)))\rangle 
+\langle \nabla _WA(x(t)),x^{(1)}(t)-A(x(t))\rangle 
~dt~=$$
$$\sum_{j=1}^{p-1}\langle W(s_j),\Delta x^{(1)}(s_j)\rangle -\int _{t_0}^{t_n}\langle W,\nabla _t x^{(1)}(t)-\frac{1}{2}{\rm grad }\Vert A\Vert ^2)\rangle -\theta _{A,x^{(1)}}(W)
~dt~=$$
$$\sum_{j=1}^{p-1}\langle W(s_j),\Delta x^{(1)}(s_j)\rangle -\int _{t_0}^{t_n}\langle W,\nabla _t x^{(1)}(t)-\frac{1}{2}{\rm grad }\Vert A\Vert ^2-\theta _{A,x^{(1)}}^{-T}\rangle 
~dt.$$
$\Box$

\spp
We say $(x,A)\in {\cal X}\times {\cal A}$ is $J${\em -critical} (or just {\em critical}) when, for all variations $\tilde x$ of $x$, 
$$\frac{\partial }{\partial h}J(x_h,A)\vert _{h=0}~=~0.$$
Then $x$ is called a {\em conditional extremum}. From Lemma \ref{lem1} follows 
\begin{theorem}\label{thm1} If $(x,A)\in {\cal X}\times {\cal A}$ is critical if and only if $x$ has singular parameters $t_0,t_1,\ldots t_n$ and, for all $k=1,2,\ldots ,n$ and all $t\in (t_{k-1},t_k)$, 
\begin{equation}\label{eleq}\nabla _t x^{(1)}(t)~=~\frac{1}{2}{\rm grad }\Vert A\Vert ^2+\theta _{A,x^{(1)}}^{-T}.\end{equation}
$\Box$
\end{theorem}

\spp
In order to be optimal, $J$ should also be critical with respect to variations in $A$. In the present paper we assume that $A$ has been found somehow; perhaps  ${\cal A}$ is a singleton. By Theorem \ref{thm1}, a conditional extremum is a $C^2$ track-sum of feasible curves satisfying (\ref{eleq}), namely conditional extrema for the case $n=1$. 

\spp
If $A={\bf 0}$, a conditional extremum is the same as a geodesic. 
\begin{example}\label{ex1} Let $n=1$ and $M=E^m$. Identifying vector fields on $E^m$ with functions $A:E^m\rightarrow E^m$, 
$$ \theta _{A,x^{(1)}}(Y) ~=~\langle \frac{d}{dt}A(x(t)),Y\rangle -\langle dA_x(Y),x^{(1)}\rangle .$$
Let $A$ be affine, of the form $A(y)=By+c$ where $c\in E^m$ and $B$ is linear. Then $x$ is critical when, for all $t\in (t_{0},t_1)$,  
$$x^{(2)}(t)~=~B^{\bf T}(Bx(t)+c) +Bx^{(1)}(t)-B^{\bf T}x^{(1)}(t)$$
where ${\bf T}$ means matrix transpose. Equivalently ~$z^{(1)}(t)+B^{\bf T}z=B^{\bf T}c$~ where ~$z(t):=x^{(1)}(t)-Bx(t)$. So   
$\displaystyle{z(t)=e^{-tB^{\bf T}}d+c}$~ where  $d\in E^m$,~ and 
$$x(t)~=~e^{(t-t_0)B}y_0+e^{tB}\int _{t_0}^te^{-sB}(e^{-sB^{\bf T}}d+c)~ds$$
where $d$ is chosen so that $x(t_1)=x_1$. Example \ref{eucex1} dealt with $B={\bf 0}$.  $\Box$
\end{example}
  
\begin{theorem}\label{cor0} Let  $x:[t_0,t_1]\rightarrow M$ be a solution of (\ref{eleq}). For some $b\in \R$ and all $t\in (t_0,t_1)$,  
\begin{equation}\label{cons}\Vert x^{(1)}(t)\Vert ^2~=~\Vert A(x(t))\Vert ^2 +b.\end{equation}
If $M$ is complete as a metric space and $A$ is $C^{\infty}$ and uniformly bounded, then $x$ extends to a unique $C^{\infty}$ solution of (\ref{eleq}) defined on all of $\R$.  
\end{theorem}

\spp
{\bf Proof:} Taking inner products with $x^{(1)}$ of both sides of (\ref{eleq}), 
$$\frac{1}{2}\frac{d}{dt}\langle x^{(1)},x^{(1)}\rangle ~=~\langle \nabla _tx^{(1)},x^{(1)}\rangle ~=~\frac{1}{2}\frac{d}{dt}\Vert A\Vert ^2+\theta _{A,x^{(1)}}(x^{(1)})~=~\frac{1}{2}\frac{d}{dt}\Vert A\Vert ^2.$$
Integrating both sides, (\ref{cons}) follows.  Alternatively, (\ref{cons}) follows from Noether's theorem and time-invariance of the Lagrangian.  

\spp
Suppose now that $M$ is complete, and that $A$ is $C^{\infty}$ and uniformly bounded. For some real $\beta $ and all $y\in M$ we have $\Vert A(y)\Vert ^2<\beta $ where $\Vert ~\Vert $ is the Riemannian norm. 
By the Picard theorem on solvability \cite{marsden} \S 7.5, $x$ is $C^{\infty}$. By (\ref{cons}),  $\Vert x^{(1)}\Vert ^2\leq \beta +b$. 
These statements hold also for any $C^2$ extension of $x$ to a solution of (\ref{eleq}).

\spp
Let ${\cal T}_+$ (respectively ${\cal T}_-$) be the nonempty sets of  
real $\tilde t_1\geq t_1$ (respectively $\tilde t_0\leq t_0$) such that $x$ extends to a $C^\infty$ solution of (\ref{eleq}) defined on $(t_0,\tilde t_1)$ (respectively $(\tilde t_0,t_1)$. To complete the proof it suffices to show that ${\cal T}_+$ is not bounded above and ${\cal T}_-$ is not bounded below. 

\spp
If ${\cal T}_+$ is bounded above, let $\bar t_1=\sup {\cal T}_+$. Then, for any integer $p\geq P$ with $P$ sufficiently large, $x$ extends to a $C^\infty $ solution $x_{[p]}:(t_0,t_1-1/p)\rightarrow M$ of (\ref{eleq}). Set $y_p:=x_{[p]}(\bar t_1-2/p)$. Because the $\Vert x_{[p]}^{(1)}\Vert $ are uniformly bounded for $p\geq P$, $\{ y_p:p\geq P\} \subset M$ is Cauchy, with limit  $\bar y$ say. Using some coordinate chart containing $\bar y$, represent each $x_{[p]}^{(1)}(\bar t_1-2/p)$ as a vector $z_p$ in $E^m$. Because the $x_{[p]}^{(1)}$ are bounded with respect to the Riemannian norm, the sequence $\{ z_p:p\geq P\} \subset E^m$ is also bounded in the Euclidean norm. So $\{ z_p:p\geq P\} $ has a convergent subsequence 
$\{ z_{p_q}:q\geq 1\} $ whose limit represents $\bar z\in TM_{\bar y}$. 

\spp
By the Picard theorem, for some $\tau >0$, any $\hat t\in \R $, and any $(\hat y,\hat z)\in TM$ sufficiently near $(\bar y,\bar z)$, 
there is a unique $C^\infty$ solution $\hat x:(\hat t-\tau ,\hat t+\tau )\rightarrow M$ of (\ref{eleq}) satisfying 
$\hat x(\hat t)=\hat y$ and $\hat x^{(1)}(\hat t)=\hat z$.  

\spp
Choose $q$ so large that, with $\hat t=\bar t_1-2/p_q$,  
\begin{itemize}
\item $(\hat y,\hat z):=(x_{[p_q]}(\hat t),x_{[p_q]}^{(1)}(\hat t))$ is sufficiently near $(\bar y, \bar z)$ and 
\item $2/p_q<\tau$. 
\end{itemize}
Splicing $x_{[p]}$ and $\hat x$ gives a $C^\infty$ extension of $x$ defined over $(t_0,\bar t_1+\tau -2/p_q)$. Since $2/p_q<\tau $ this contradicts the definition of $\bar t_1$ as $\sup {\cal T}_+$. The proof that ${\cal T}_-$ is not bounded below is entirely similar. $\Box$

\spp
\begin{example}\label{Abdex} For $A={\bf 0}$ the hypothesis that $M$ be complete is needed for extendability, from the Hopf-Rinow Theorem. Extendability may also fail when $M$ is complete and $A$ is unbounded, as when   $M=E^1$ and $A(x)=x$. $\Box$
\end{example}
\section{Closed Prior Fields}\label{closedsec}
The definition of $J(x,A)$ depends not only on $x\in {\cal X}$ and $A\in {\cal A}$, but also on the $t_k\in \R$ and the $x_k\in M$ which are used to define ${\cal X}$. 
The {\em reverse} $\bar x$ of $x$ is defined by $\bar x(u):=x(-u)$ for $u\in [-t_n,-t_0]$ and 
$$\bar {\cal X}~:=~{\cal X}_{-t_n,-t_{n-1},\ldots ,-t_0;x_n,x_{n-1},\ldots x_0},\quad  
\bar J~:=~J_{-t_n,-t_{n-1},\ldots ,-t_0;x_n,x_{n-1},\ldots x_0}:\bar {\cal X}\times {\cal A}\rightarrow \R .$$
 If $x\in {\cal X}$ then $\bar x\in \bar {\cal X}$. From Theorem \ref{thm1} follows 
 \begin{corollary}\label{prop1} Suppose $dA^T={\bf 0}$. Then $(x,A)\in {\cal X}\times {\cal A}$ is $J$-critical if and only if $(\bar x,A)\in \bar {\cal X}\times {\cal A}$ is $\bar J$-critical.  
 \end{corollary}
 
 \spp
 {\bf Proof:} $dA^T={\bf 0}\Longrightarrow \theta _{A,x^{(1)}}={\bf 0}=\theta _{A,\bar x^{(1)}}$, and  so $(x,A)$ is critical if and only if 
 $$\nabla _tx^{(1)}(t)~=~\frac{1}{2}{\rm grad}(\Vert A\Vert ^2)~\Longleftrightarrow ~\nabla _{\bar t}\bar x^{(1)}(\bar t)~=~\frac{1}{2}{\rm grad}(\Vert A\Vert ^2).$$
 $\Box$
 
\spp
However if $dA^T={\bf 0}$ and $(x,A)$ is $J$-optimal, $(\bar x ,A)$ need not be $\bar J$-optimal. 
\begin{example}\label{counterprop1} Let $A$ be the clockwise unit vector field on $M=S^1$ with the standard Riemannian metric. Take $n=1$, $t_0=0$, $t_1=2\pi $ and $x_0=x_1=(1,0)$. Then $(x,A)$ given by $x(t)=(\cos t,-\sin t)$ is $J$-optimal, with $J(x,A)=0$, and $\bar J(\bar x,A)=8\pi $. Also $(x^*,A)$ given by $x^*(u)=(\cos u,-\sin u)$ for $u\in [-2\pi ,0]$ is $\bar J$-optimal, with $\bar J(x^*,A)=0$. $\Box$ 
\end{example}
 %
 %
 %
 \begin{definition}\label{def1}
 A nonempty set ${\cal A}$ of $C^1$ vector fields on $M$ is said to be {\em reflexive} when, for any real $t_0<t_1<\ldots  <t_n$, and any $x_0,x_1,\ldots ,x_n\in M$, 
there exists $c_{A,x_0,x_n}\in \R$ independent of $t_0,t_1,\ldots ,t_n$ and independent of $x_1,x_2,\ldots ,x_{n-1}$, and there exists $\bar A\in {\cal A}$ depending only on $A$, such that, for all $x\in {\cal X}$,
$\bar J(\bar x,\bar A)=J(x,A)+c_A$. 
Call $A$ {\em reflexive} when $\{ A\} $ is reflexive. $\Box$
\end{definition} 
In particular, if ${\cal A}$ is reflexive then for all $(x,A)\in {\cal X}\times {\cal A}$ 
$$(x,A)\in {\cal X}\times {\cal A}\hbox{~is~}J\hbox{-critical}\quad \Longleftrightarrow \quad  (\bar x,\bar A)\in \bar {\cal X}\times {\cal A}\hbox{~is~}\bar J\hbox{-critical}$$
and
$$(x,A)\in {\cal X}\times {\cal A}\hbox{~is~}J\hbox{-optimal}\quad \Longleftrightarrow \quad  (\bar x,\bar A)\in \bar {\cal X}\times {\cal A}\hbox{~is~}\bar J\hbox{-optimal}. $$
If, as in Example \ref{geoex}, ${\cal A}$ is closed under multiplication by $-1$ then ${\cal A}$ is reflexive, with $c_{A,x_0,x_n}=0$ and $\bar A=-A$. In Example \ref{eucex1} any nonempty set of constant vector fields on $E^m$ is reflexive, with $c_{A,x_0,x_n}=4\langle A,x_n-x_0\rangle $ and $\bar A=A$. The condition that ${\cal A}$ be reflexive is especially stringent when ${\cal A}$ is a singleton.  
\begin{example}\label{refex} For $M$ the unit $2$-sphere $S^2$ in $E^3$, set ${\cal A} =\{ A\}$ 
where $A(v_1,v_2,v_3)=(v_2,-v_1,0)$ for $v\in S^2$. The integral curves of $A$ are shown (black) as latitudinal circles traversed in the clockwise direction. Set $n=1$, $t_0=0$, $t_1=1$, $x_0=(0.866,0,0.5)$ and $x_1=(0.5187,0.8486,0.1039)$. Figure \ref{figlat} 
 shows (red) the conditional minimum $x$ and (blue) the conditional minimum $x^*$ for the reverse data, as approximated from a numerical computation. Theorem \ref{horthm} gives closed form expressions for such curves. 
 We find $\bar J(x^*,A)\approx 0.18<\bar J(\bar x,A)\approx 0.44$. So $A$ is not reflexive. 
\begin{figure}[htbp] 
   \centering
   \includegraphics[width=2in]{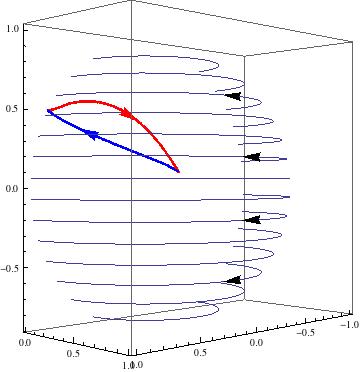} 
   \caption{Conditional Minima With $A$ Not Reflexive, in Example \ref{refex}}
   \label{figlat}
\end{figure}
$\Box$
\end{example}

\spp
Indeed in Example \ref{refex} $A$ cannot be reflexive, by the first part of
\begin{theorem}\label{thm2} If $A$ is reflexive then $dA^T={\bf 0}$. If $A^T=d\phi $ where $\phi :M\rightarrow \R $ then $A$ is reflexive, with $c_{A,x_0,x_n}=4(\phi (x_n)-\phi (x_0))$. 
\end{theorem}

\spp
{\bf Proof:} Suppose first that $A$ is reflexive, and choose $t_0<t_1<\ldots <t_n$ and $x_0,x_1,\ldots ,x_n\in M$. For any $x\in {\cal X}$, and any variation $\tilde x$ of $x$, $2K(h):=J(x_h,A)-\bar J(\bar x_h,A)=c_{A,x_0,x_n}$. By Lemma \ref{lem1},~~ 
$0=K'(0)=$ 
$$\sum_{j=1}^{p-1}\langle W(s_j),\Delta x^{(1)}(s_j)\rangle -\int _{t_0}^{t_n}\langle W,\nabla _t x^{(1)}(t)-\frac{1}{2}{\rm grad }\Vert A\Vert ^2-\theta _{A,x^{(1)}}^{-T}\rangle 
~dt~$$
$$-\sum_{j=1}^{p-1}\langle \bar W(-s_j),\Delta \bar x^{(1)}(-s_j)\rangle +\int _{-t_n}^{-t_0}\langle \bar W(\bar t),\nabla _{\bar t} \bar x^{(1)}(\bar t)-\frac{1}{2}{\rm grad }\Vert A\Vert ^2-\theta _{A,\bar x^{(1)}}^{-T}\rangle 
~d\bar t~=~$$
$$-\int _{t_0}^{t_n}\langle W(t),\nabla _t x^{(1)}(t)-\frac{1}{2}{\rm grad }\Vert A\Vert ^2-\theta _{A,x^{(1)}}^{T}\rangle 
~dt~+~\int _{t_0}^{t_n}\langle W(t),\nabla _{t} x^{(1)}(t)-\frac{1}{2}{\rm grad }\Vert A\Vert ^2+\theta _{A,x^{(1)}}^{-T}\rangle 
~dt~=~$$
$$2\int _{t_0}^{t_n}\langle W(t),\theta _{A,x^{(1)}}^{-T}\rangle ~dt.$$
Since this holds for all variations $\tilde x$,  we have $\theta _{A,x^{(1)}}\vert _{x(t)}={\bf 0}$ for all nonsingular $t$ and, by continuity, for all $t\in [t_0,t_n]$. Since this holds for all $x\in {\cal X}$, $dA^T={\bf 0}$. \\

\spp
Suppose now that $\omega _A=d\phi $. Then $J(x,A)-\bar J(\bar x,A)=$
$$\int _{t_0}^{t_n}\Vert x^{(1)}(t)-A(x(t))\Vert ^2~-~\Vert x^{(1)}(t)+A(x(t))\Vert ^2~dt~=~-4\int_{t_0}^{t_n} \langle A(x(t)),x^{(1)}(t)\rangle ~dt~=$$
$$-4\int _{t_0}^{t_n}\frac{d}{dt}(\phi \circ x(t))~dt~=~-4(\phi (x_n)-\phi (x_0)).$$
$\Box$

\spp
\begin{corollary}\label{cohomcor} Let $H^1(M;\R )={\bf 0}$. Then $A$ is reflexive if and only if $A$ is conservative. $\Box$
\end{corollary}
\begin{example} For $p\geq 2$ let $M$ be any Riemannian manifold with the homotopy type of the unit sphere $S^p$ in $E^{p+1}$, or real projective $p$-space $\R P^p$, or complex projective $m$-space $\C P^p$. Then $A$ is reflexive if and only if $A$ is conservative. In particular, regardless of the Riemannian metric on  $M=SO(3)$, $A$ is reflexive if and only if $A$ is conservative. $\Box$
\end{example}
%
\section{Left-Invariant Priors on Bi-Invariant Lie Groups}\label{Liesec}  
Take  $M$ to be a Lie group $G$ with a bi-invariant Riemannian metric $\langle ~,~\rangle $, and let $A$ be a left-invariant vector field on $G$ \cite{milnor}, \cite{holm}, \cite{vara}. Let $n=1$. The {\em left Lie-reduction} of a vector field $Z$ defined along $x:[t_0,t_1]\rightarrow G$ is defined to be the curve $Z_L$ in the Lie algebra ${\cal G}=TG_{\bf 1}$ given by 
$$Z(x(t))~=:~dL(x(t))_{\bf 1}(Z_L(t))$$
where $L$ is left-multiplication, and ${\bf 1}$ is the identity element  
of $G$. Then $A_L$ is constant. Denote $(x^{(1)})_L$ by $V$. 
Because $\langle ~,~\rangle $ and $A$ are left-invariant,  for any $Y\in TG_{x(t)}$,  
$$\langle Y_L,(\theta _{A,x^{(1)}}^{-T})_L\rangle ~=~\langle Y,\theta _{A,x^{(1)}}^{-T}\rangle ~=
~\theta _{A,x^{(1)}}(Y)~=~
(dA^T)(x^{(1)},Y)~=~\langle \nabla _tA,Y\rangle -\langle \nabla _YA,x^{(1)}\rangle ~=~$$
$$\langle (\nabla _tA)_L,Y_L\rangle -\langle (\nabla _YA)_L,V\rangle ~=~
\frac{1}{2}\langle [V,A_L],Y_L\rangle -\frac{1}{2}\langle [Y_L,A_L],V\rangle ~=~\langle [V,A_L],Y_L\rangle $$
namely $\displaystyle{(\theta _{A,x^{(1)}}^{-T})_L=[V,A_L]}$. 
Take $n=1$. Then a feasible curve $x$ is a conditional extremal when $\displaystyle{V^{(1)}(t)=[V(t),A_L]}$. Equivalently  ~$\displaystyle{V(t)={\rm Ad}(e^{-(t-t_0)A_L})V(t_0)}$, and 
$$\frac{d}{dt}(x(t)e^{-(t-t_0)A_L})~=~dR(e^{-(t-t_0)A_L})_{x(t)}\circ dL(x(t))_{\bf 1}V(t)-dL(x(t)e^{-(t-t_0)A_L})_{\bf 1}A_L~=~$$
$$~dR(e^{-(t-t_0)A_L})_{x(t)}\circ dL(x(t))_{\bf 1}\circ dL(e^{-(t-t_0)A_L})_{e^{(t-t_0)A_L}}\circ dR(e^{(t-t_0)A_L})_{\bf 1}V(t_0)-dL(x(t)e^{-(t-t_0)A_L})_{\bf 1}A_L$$
$$~=~dL(x(t)e^{-(t-t_0)A_L})_{\bf 1}(V(t_0)-A_L).$$
So $x(t)e^{-(t-t_0)A_L}=e^{(t-t_0)(V(t_0)-A_L)}$. Taking $B_L:=V(t_0)-A_L$, this proves the first assertion of 
\begin{theorem}\label{thm2g} Let $A$ be a left-invariant vector field on a Lie group $G$ with a bi-invariant Riemannian metric.  Let $n=1$. A feasible curve $x:[t_0,t_1]\rightarrow G$ is a conditional extremum 
if and only if, for all $t\in [t_{0},t_1]$, 
\begin{equation}\label{xeq}x(t)~=~e^{(t-t_{0})B_L}e^{(t-t_{0})A_L}x_{0}\end{equation}
where $B_L\in {\cal G}$ is such that $x(t_1)=x_1$. Then, for all $t\in (t_0,t_1)$, 
$$\Vert x^{(1)}(t)\Vert ~=~\Vert A_L+B_L\Vert ,\quad \Vert x^{(1)}(t)-A(x(t))\Vert~=~\Vert B_L\Vert ,\quad   
\langle x^{(1)}(t),A(x(t))\rangle ~=~\langle A_L+B_L,A_L\rangle .$$ 
In particular 
$\displaystyle{J(x,A)=(t_1-t_{0})\Vert B_L\Vert ^2}$.  
\end{theorem}

\spp
{\bf Proof:} Let $x$ be a conditional extremum. Because $\langle ~,~\rangle $ is left-invariant, 
$\Vert x^{(1)}(t)\Vert =\Vert V(t)\Vert =\Vert V(t_0)\Vert $, 
since ${\rm Ad}$ acts by isometries. So $\Vert x^{(1)}(t)\Vert =\Vert A_L+B_L\Vert $. 
Because $A$ is left-invariant and $V(t)={\rm Ad}(e^{-(t-t_{0})A_L})V(t_{0})$,
$$\Vert V(t)-A_L\Vert ~=~ 
\Vert {\rm Ad}(e^{-(t-t_{0})A_L})V(t_{0})-A_L\Vert ~=~\Vert V(t_{0})-A_L\Vert ~=~\Vert B_L\Vert $$
because the isometry ${\rm Ad}(e^{(t-t_{0})A_L}):{\cal G}\rightarrow {\cal G}$ fixes $A_L$.   
Similarly, $\langle x^{(1)}(t),A(x(t))\rangle $ is 
$$~\langle V(t),A_L\rangle ~=~\langle {\rm Ad}(e^{-(t-t_{0})A_L})V(t_{0}),A_L\rangle =~
\langle V(t_0),A_L\rangle ~=~\langle A_L+B_L,A_L\rangle .$$
$\Box$
\begin{example} Suppose $x_1=e^{t_*A_L}x_0$ for some $t_*>t_0$. Substituting $\displaystyle{B_L=\frac{t_*-t_1}{t_1-t_0}A_L}$ in (\ref{xeq}), we find that the geodesic given by 
$$\displaystyle{x(t)~=e~^{\frac{(t-t_0)(t_*-t_0)}{t_1-t_0}A_L}x_0}$$
is a conditional extremum. $\Box$
\end{example}
\begin{corollary} Let $A$ be a left-invariant vector field on a bi-invariant Lie group $G$. Then $x:\R \rightarrow G$ satisfies equation (\ref{eleq}) with $x(0)={\bf 1}$ if and only if $x$ is a pointwise product 
$$x(t)~=~e^{tB_L}e^{tA_L}$$
of one-parameter subgroups of $G$. $\Box$  
\end{corollary}
\begin{example}\label{S3ex} Take $G$ to be the group $S^3\subset \H $ of unit quaternions, with bi-invariant Riemannian metric from $E^4\cong \H $.  A left-invariant vector field $A$ corresponds to a pure imaginary $\alpha \in \H$, where $A(x)=x\alpha $. The 
one-parameter subgroups of $S^3$ are given by $g_{\beta }(t)=e^{t{\beta }}$ where ${\beta }\in \H$ is pure imaginary. 
Choose $q_\alpha ,q_\beta \in S^3$ so that $q_\alpha \alpha \bar q_\alpha =a\imath $
and $q_\beta \beta \bar q_\beta =b\imath $ where $a,b\in \R $ and $\bar q$ is the conjugate of $q\in \H$.  Solutions $x:\R \rightarrow S^3$ of (\ref{eleq}) with $x(0)={\bf 1}$ have the form 
$$x(t)~=~\bar q_\alpha (\cos (at)+\imath \sin (at))q_\alpha \bar q_\beta        (\cos (bt)+\imath \sin (bt))q_\beta .$$
Since $A$ is nowhere-zero and $M$ is compact, $A$ is not conservative.   

\spp
Take $A_L=(-0.5,-0.5,0.3)$, $x_0={\bf 1}$ and $x_1=(-0.0359448, -0.228089, -0.937324, -0.260972)$. 
Figure \ref{s31fig} shows (red) ${\rm exp}^{-1}\circ x:[0,\pi ]\rightarrow E^3$ for a conditional minimum $x:[0,\pi ]\rightarrow S^3$.  The preimage under ${\rm exp}$ of the orbit $y:[0,\pi ]\rightarrow S^3$ of $A$ is shown (black) as a directed line segment, with (blue) the preimage of $t\mapsto e^{tB_L}y(\pi )$ for $t\in [0,\pi ]$, where $B_L=(0.2,0.2,0.2)$. 
\begin{figure}[htbp] 
   \centering
   \includegraphics[width=2in]{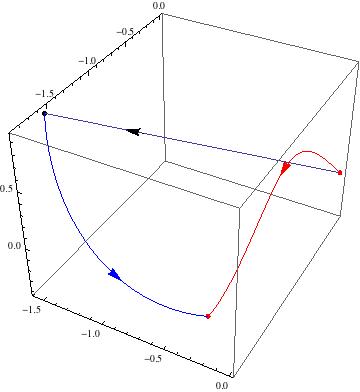} 
   \caption{A Pointwise Product of $1$-Parameter Subgroups of $S^3$, in Example \ref{S3ex}}
   \label{s31fig}
\end{figure}
$\Box$
\end{example}
For $1\leq k\leq n$ set $s_k:=t_k-t_{k-1}$ and $y_k:=x_kx_{k-1}^{-1}$. 
\begin{corollary}\label{cor11} Let $(x,A)$ be optimal, where $A$ is a left-invariant vector field on a Lie group $G$ with bi-invariant Riemannian metric.  For some $B_{L,1},B_{L,2},\ldots ,B_{L,n}\in {\cal G}$,
$$x(t)~=~e^{(t-t_{k-1})B_{L,k}}~e^{(t-t_{k-1})A_L}x_{k-1}\quad \hbox{for ~}t\in [t_{k-1},t_k]$$ 
where $e^{s_kB_{L,k}}e^{s_kA_L}=y_k$, and the geodesic arcs  
$\displaystyle{s\in [0,s_k]\mapsto e^{sB_{L,k}}\in G}$ are minimal. 
$\Box$
\end{corollary}

\spp
When ${\cal A}$ is finite, Corollary \ref{cor11} allows us to compute all optimal pairs $(x,A)$ in terms of the exponential map of $G$. In particular, if  for each $k$ there is a unique minimal geodesic on $[0,s_k]$ joining $e^{s_kA_L}$ to $y_k$, then $x$ is uniquely defined by $A$.  
%

\subsection{Manifolds of Left-Invariant Priors}\label{lsr}
Let ${\cal A}$ be a set of left-invariant vector fields on $G$ where ${\cal A}_L:=\{ A_L:A\in {\cal A}\}$ is a $C^1$ immersed submanifold of ${\cal G}$ of dimension $p\leq m$. The following analysis can be developed for when ${\cal A}_L$ is the image of a $C^1$ immersion (not necessarily one-to-one), but we are interested in cases where ${\cal A}_L$ is an affine subspace of ${\cal G}$, especially ${\cal A}_L={\cal G}$. 
\begin{theorem}\label{sumBthm} Let $(x,A)\in {\cal X}\times {\cal A}$ be optimal, where $e^{-s_kA_L}y_k$ is not a critical value of the exponential map for any $k$. Defining $B_{L,k}$ in terms of $A_L$, as in Corollary \ref{cor11}, 
$$\bar B_L~:=~\sum_{k=1}^ns_k\left(  \frac{{\bf 1}-e^{-s_k{\rm ad}A_L}}{{\rm ad}A_L} \right)B_{L,k}$$
is orthogonal to the tangent space of ${\cal A}_L$ at $A_L$. 
\end{theorem}

\spp
{\bf Proof:} For $\tilde A_L$ near $A_L$, the relation 
\begin{equation}\label{tilrel}e^{s_k\tilde B_{L,k}}~=~y_k~e^{-s_k\tilde A_L}~\end{equation}
determines $\tilde B_{L,k}$ locally as a $C^\infty$ function of $\tilde A_L$, with value $B_{L,k}$ at $A_L$.  
Using the formula for the derivative of the exponential map \cite{tuyn}, \cite{vara} to differentiate (\ref{tilrel}) at $A_L$, 
$$dL(e^{s_kB_{L,k}})\circ \left( \frac{{\bf 1}-e^{-s_k{\rm ad}B_{L,k}}}{{\rm ad}B_{L,k}}  \right) \circ d\tilde B_{L,k}(W)~=~dL(y_k)\circ dL(e^{-s_kA_L})\circ \left(  \frac{{\bf 1}-e^{s_k{\rm ad}A_L}}{{\rm ad}A_L} \right) (W)   $$
where $L(g)$ is left multiplication by $g\in G$, and $W$ is tangent to ${\cal A}_L$ at $A_L$. Then by (\ref{tilrel}), 
$$\left( \frac{{\bf 1}-e^{-s_k{\rm ad}B_{L,k}}}{{\rm ad}B_{L,k}}  \right) \circ d\tilde B_{L,k}(W)~=~\left(  \frac{{\bf 1}-e^{s_k{\rm ad}A_L}}{{\rm ad}A_L} \right) (W)   ~\Longrightarrow $$
$$\langle B_{L,k},d\tilde B_{L,k}(W)\rangle ~=~\langle B_{L,k},\left( \frac{{\bf 1}-e^{-s_k{\rm ad}B_{L,k}}}{{\rm ad}B_{L,k}}  \right) \circ d\tilde B_{L,k}(W)\rangle~=$$
$$~\langle B_{L,k},\left(  \frac{{\bf 1}-e^{s_k{\rm ad}A_L}}{{\rm ad}A_L} \right) (W)\rangle  ~=~
~-\langle \left(  \frac{{\bf 1}-e^{-s_k{\rm ad}A_L}}{{\rm ad}A_L} \right)B_{L,k}, W\rangle .$$
By Theorem \ref{thm2g}, $A_L$ is a critical point of $J(x,A)=\sum_{k=1}^ns_k\Vert \tilde B_{L,k}(A_L)\Vert ^2$. So 
$$\langle \sum_{k=1}^ns_k\left(  \frac{{\bf 1}-e^{-s_k{\rm ad}A_L}}{{\rm ad}A_L} \right)B_{L,k}, W\rangle ~=~0$$
for all $W$ tangent to ${\cal A}_L$ at $A_L$.  $\Box$

\spp
\begin{corollary} Let ${\cal A}_L={\cal G}$. Then $\displaystyle{\langle A_L,\sum_{k=1}^ns_kB_{L,k}\rangle =~0}$. 
\end{corollary}

\spp
{\bf Proof:} Since $\bar B_L={\bf 0}$,~ 
$\displaystyle{0=\langle A_L,\sum_{k=1}^ns_k\left(  \frac{{\bf 1}-e^{-s_k{\rm ad}A_L}}{{\rm ad}A_L} \right)B_{L,k}\rangle =}$
$$-\sum_{k=1}^ns_k\langle \left(  \frac{{\bf 1}-e^{s_k{\rm ad}A_L}}{{\rm ad}A_L} \right)A_L,B_{L,k}\rangle ~=~
-\sum_{k=1}^ns_k\langle A_L,B_{L,k}\rangle .$$
$\Box$
\begin{corollary}\label{corsum} Let ${\cal A}_L={\cal G}$, $s_1=s_2=\ldots =s_n$ and suppose $s_1A_L$ is not a critical point of the exponential map. Then, with the hypotheses of Theorem \ref{sumBthm},
$$\sum_{k=1}^nB_{L,k}~=~{\bf 0}.$$ 
\end{corollary}

\spp
{\bf Proof:} The derivative of the exponential map at $s_1A_L$ is $\displaystyle{dL(e^{s_1A_L})\circ \left( \frac{{\bf 1}-e^{-s_1{\rm ad}A_L}}{s_1{\rm ad}}\right) }$. $\Box$
\begin{example}\label{S3ex2} Let $G=S^3$, ${\cal A}_L={\cal G}$, $n=4$ and $s_1=s_2=s_3=s_4=1/4$. Figure \ref{s3ex2fig} shows (red) the preimage under ${\rm exp}$ of 
$x:[0,1]\rightarrow S^3$ with observations $x_0=(1,0,0,0)$, 
\begin{eqnarray*}
x_1&=&(0.1304,0.7923,0.4574,0.3821)\\
x_2&=&(0.5809, 0.0381, 0.3385, 0.7393)\\
x_3&=&(0.5523, 0.6251, 0.5513,0.0172)\\
x_4&=&(0.2810, 0.1241, 0.6817,0.6640)
\end{eqnarray*}
made at $t_0=0$, $t_1$, $t_2$, 
$t_3$ and $t_4$. Using Mathematica's {\em FindMinimum} on $\sum_{k=1}^4\Vert \tilde B_{L,k}\Vert ^2$, considered as a function of $\tilde A$, the best estimating geodesic from $x_0$ is found to have infinitesimal generator  $A_L=(1.40398, 0.196766, 1.05334)$, with preimage shown in black in Figure \ref{s3ex2fig}. 
\begin{figure}[htbp] 
   \centering
   \includegraphics[width=2in]{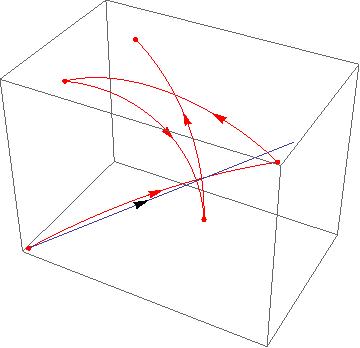} 
   \caption{The Conditional Minimum and Optimal Prior in Example \ref{S3ex2}}
   \label{s3ex2fig}
\end{figure}
We find 
\begin{eqnarray*}
B_{L,1}&=&~~(2.7669, 2.3129, 2.0736)\\
B_{L,2}&=&-(1.0075, 4.2867, 1.4298)\\
B_{L,3}&=&-(2.1680, -3.6097, 2.4150)\\
B_{L,4}&=&~~(0.4086, -1.6359, 1.7713)
\end{eqnarray*}
with sum ${\bf 0}$ in accordance with Corollary \ref{corsum}. When, as presently, observations tend to contradict the hypothesis that $x$ is an integral curve of a field near $A$, the usefulness of the optimal interpolant is questionable. For less contradictory observations we find  the preimage of $x$ in $E^3$ more nearly piecewise-affine, consistent with Example \ref{eucex1}. 
$\Box$
\end{example}
\begin{corollary}\label{3pts} Take $n=2$ in Corollary \ref{corsum}. If $(x,A)$ is optimal then 
$e^{s_1A_L}$ is the midpoint of a geodesic  joining $y_1$ and $y_2$.  According as $t\in [t_0,t_1]$ or $t\in [t_1,t_2]$,
\begin{eqnarray*}
x(t)&=&e^{(t-t_0)B_{L,1}}~e^{(t-t_0)A_L}~x_0\quad \hbox{or}\\
x(t)&=&e^{(t_1-t)B_{L,1}}~e^{(t-t_1)A_L}~x_1
\end{eqnarray*}
where $\displaystyle{e^{s_1B_{L,1}}=y_1e^{s_1A_L}}$. 
\end{corollary}

\spp
{\bf Proof:} By Corollary \ref{cor11}, for an optimal field $A$, $e^{s_1B_{L,1}}e^{s_1A_L}=y_1$ and $e^{-s_1B_{L,1}}e^{s_1A_L}=y_2$. Then $s\in [-s_1,s_1]\mapsto  e^{sB_{L,1}}e^{s_1A_L}\in G$ is a geodesic from $y_2$ to $y_1$. $\Box$

\section{Symmetric Priors on Sphere and Hyperboloid}\label{sphsec}
Define a nondegenerate symmetric bilinear form $\langle ~,~\rangle $ on $\R ^3$ by 
$$\langle (v_1,v_2,v_3),(w_1,w_2,w_3)\rangle =v_1w_1+v_2w_2+\sigma v_3w_3 $$
where $\sigma =\pm 1$. According as $\sigma$ is $1$ or $-1$, $\langle ~,~\rangle$ is the Euclidean metric or the Lorentz metric. In either case $\langle ~,~\rangle $ restricts to a Riemannian metric of constant curvature $\sigma $ on  
$$M_\sigma ~=~\{ y\in \R ^3:\langle y,y\rangle =\sigma \} .$$
Then $M_+$ is the unit $2$-sphere $S^2$ with the standard Riemannian metric, and $M_-$ is the $2$-sheeted unit hyperboloid $H^2$ \cite{ratty}. 

\spp
For vector fields $X,Y$ on $M_\sigma $ and a $C^1$ vector field $A$ on $M_\sigma $,
$$\theta _{A,X}(Y)~=~\langle \nabla _XA,Y\rangle -\langle \nabla _YA,X\rangle ~=~~\langle X(A),Y\rangle -\langle Y(A),X\rangle .$$
Define vector fields $B$ and $C$ on $M_\sigma $ by 
$B(x):=(-x_2,x_1,0)$ and $C(x):=(0,0,1)-x_3x$.  
Since $C$ is the gradient of the height function, $\theta _{C,X}={\bf 0}$.  
Also   
$$\theta _{B,X}(Y)~=~2(X_1Y_2-X_2Y_1)~=~2\langle  (-X_2,X_1,0)-\sigma (-x_1X_2+x_2X_1)x, Y\rangle .$$ 
So ~$\displaystyle{\theta _{B,X}^{-T}=2((-X_2,X_1,0)-\sigma (-x_1X_2+x_2X_1)x)}$. Given $C^1$ functions $\beta ,\gamma :M_\sigma \rightarrow \R $, set 
$$A~:=~\beta B+\gamma C~~\Longrightarrow ~~\Vert A\Vert ^2~=~(\beta ^2+\gamma ^2)(x_1^2+x_2^2).$$
 By (\ref{prodth}), 
$\displaystyle{\theta _{A,X}~=~\beta \theta _{B,X}+X(\beta )B^T-\langle B,X\rangle d\beta  +
X(\gamma )C^T-\langle C,X\rangle d\gamma }$. 

\spp
Let $A$ be invariant with respect to rotations of the first two coordinates of $M_\sigma $, namely $\beta =\bar \beta \circ pr_3$ and $\gamma =\bar \gamma \circ pr_3$ where $\bar \beta , \bar \gamma :\R \rightarrow \R $ and $pr_3:M_\sigma \rightarrow \R$ projects to the third coordinate.  Then  
$$\frac{1}{2}{\rm grad}(\Vert A\Vert ^2)=((\bar \beta ^2+\bar \gamma ^2)x_1,(\bar \beta ^2+\bar \gamma ^2)x_2,(\bar \beta \bar \beta '+\bar \gamma \bar \gamma ')(1-x_3^2))
-(\bar \beta ^2+\bar \gamma ^2+(\bar \beta \bar \beta '+\bar \gamma \bar \gamma ')x_3)(1-x_3^2)x.$$
$$\hbox{Also}\quad 
\theta _{A,X}^{-T}~=~\beta \theta _{B,X}^{-T}+X(\beta )B-\langle B,X\rangle (d\beta  )^{-T}~=~$$
$$2\bar \beta ((-X_2,X_1,0)+wx)~+
\bar \beta '((-x_2X_3,x_1X_3,-w)+wx_3x)$$
where $w:=\sigma (x_1X_2-x_2X_1)$. 
Taking $X=x^{(1)}$ for a $C^2$ curve $x:\R \rightarrow M_\sigma $, equation (\ref{eleq}) is equivalent to 
\begin{eqnarray}
\nonumber x_1^{(2)}+\sigma \Vert x^{(1)}\Vert ^2x_1&=&~~(\bar \beta ^2+\bar \gamma ^2)x_3^2x_1-
(\bar \beta \bar \beta '+\bar \gamma \bar \gamma ')(1-x_3^2)x_3x_1+\\
\label{x1eq}&~&-2\bar \beta x_2^{(1)}- 
\bar \beta 'x_2x_3^{(1)} +(2\bar \beta +\bar \beta 'x_3)wx_1 \\
\nonumber x_2^{(2)}+\sigma \Vert x^{(1)}\Vert ^2x_2&=&~~(\bar \beta ^2+\bar \gamma ^2)x_3^2x_2-
(\bar \beta \bar \beta '+\bar \gamma \bar \gamma ')(1-x_3^2)x_3x_2+ \\
\label{x2eq}&~& 2\bar \beta x_1^{(1)}+ \bar \beta 'x_1x_3^{(1)} +(2\bar \beta + \bar \beta 'x_3)wx_2 \\
\nonumber x_3^{(2)}+\sigma \Vert x^{(1)}\Vert ^2x_3&=& -(\bar \beta ^2+\bar \gamma ^2)(1-x_3^2)x_3+(\bar \beta \bar \beta '+\bar \gamma \bar \gamma ')(1-x_3^2)^2+\\
\label{x3eq}&~&(2\bar \beta x_3-\bar \beta '(1-x_3^2))w.
\end{eqnarray}
By equation (\ref{cons}), for some $b\in \R$, 
\begin{equation}\label{conssph}
\sigma \Vert x^{(1)}\Vert ^2~=~(\bar \beta ^2+\bar \gamma ^2)(1-x_3^2)+b.
\end{equation}
By rotational symmetry and Noether's theorem,  
\begin{equation}\label{consrot}
w~=~\sigma (x_1x_2^{(1)}-x_2x_1^{(1)})~=~\bar \beta (1-x_3^2)+c.
\end{equation}
where $c\in \R $. Set $r(t)=\sqrt{x_1^2+x_2^2}=\sqrt{\sigma (1-x_3^2)}$. For $r\vert (a_0,a_1)$ nowhere-zero, choose $\psi :(a_0,a_1)\rightarrow \R$ so that 
$x_1=r\cos \psi$ and $x_2=r\sin \psi $. 
Then $w=\sigma r^2\psi ^{(1)}$ and, by (\ref{consrot}), 
\begin{equation}\label{psieq}
\psi ^{(1)}(t)~=~\bar \beta +\frac{c}{1-x_3^2}.
\end{equation}
So $x_3:(a_0,a_1)\rightarrow \R$ determines $x\vert (a_0,a_1)$ up to quadrature. 
\begin{proposition}\label{constintegrand} $\displaystyle{\sigma \Vert x^{(1)}(t)-A(x(t))\Vert ^2=b-2c\bar \beta 
+2\bar \gamma (1-2x_3^2)x_3^{(1)}}$.
\end{proposition}

\spp
{\bf Proof:} $\displaystyle{\Vert x^{(1)}(t)-A(x(t))\Vert ^2=\Vert x^{(1)}-\bar \beta (-x_2,x_1,0)
-\bar \gamma (-x_3x_1,-x_3x_2,1-x_3^2)\Vert ^2=}$
$$\Vert x^{(1)}\Vert ^2+\sigma (\bar \beta ^2+\bar \gamma ^2)(1-x_3^2)-2\sigma \bar \beta w
+2\bar \gamma (x_3x_1x_1^{(1)}+x_3x_2x_2^{(1)}+\sigma (1-x_3^2)x_3^{(1)})$$
$$~=~\sigma b-2\sigma \bar \beta c
+2\bar \gamma (-\sigma x_3^2x_3^{(1)}+\sigma (1-x_3^2)x_3^{(1)})~=~\sigma b-2\sigma c\bar \beta 
+2\sigma \bar \gamma (1-2x_3^2)x_3^{(1)}$$
by (\ref{conssph}), (\ref{consrot}), and because $x_1^{(1)}x_1+x_2^{(1)}x_2+\sigma x_3^{(1)}x_3=0$. $\Box$

\begin{corollary} Let $\bar \beta $ and $\bar \gamma $ be constant. Then 
$$\sigma \int \Vert x^{(1)}(t)-A(x(t))\Vert ^2~dt~=~(b-2\bar \beta c)t
+2\bar \gamma (x_3-2x_3^3/3)+c_0$$
where $c_0$ is constant. $\Box$
\end{corollary}
\begin{proposition}\label{corconstx3} Let $\bar \beta $ and $\bar \gamma $ be constant. A solution $x:\R \rightarrow M_\sigma $ of (\ref{x1eq}), (\ref{x2eq}), (\ref{x3eq}), $x_3$ is constant if and only if, for some $\psi _0\in \R$, either 
\begin{enumerate}
\item $x(t)=(\cos (\omega t+\psi _0),\sin (\omega t+\psi _0),0)$ for all $t$, where $\omega \in \R$ is arbitrary, or 
\item $\bar \gamma =0$ and $x(t)=(\sqrt{\sigma (1-h^2)}\cos (\bar \beta t+\psi _0),\sqrt{\sigma (1-h^2)}\sin (\bar \beta t+\psi _0),h)$ 
for any $h$ with $\sigma (1-h^2)>0$.
\end{enumerate}
\end{proposition}

\spp
{\bf Proof:} Suppose $x_3(t)=h$ for all $t$. If $h^2= 1$ then $x$ is constant at $\pm (0,0,1)$. Conversely these constant solutions satisfy (\ref{x1eq}), (\ref{x2eq}), (\ref{x3eq}). Also $b=c=0$. 
Suppose now $h^2\not= 1$ and $\bar \beta \not= 0$. By (\ref{psieq}), 
$$x(t)~=~(\sqrt{\sigma (1-h^2)}\cos (\omega t+\psi _0),\sqrt{\sigma (1-h^2)}\sin (\omega t+\psi _0),h)$$
where $\displaystyle{\omega :=\bar \beta +\frac{c}{1-h^2}}$ and $\psi _0\in \R $. Equations (\ref{x1eq}), (\ref{x2eq}) are together equivalent to  
$$c^2+2(1-h^2)^2\bar \beta c-((1-2h^2)\bar \gamma ^2+b)(1-h^2)^2~=~0,$$
and equation (\ref{x3eq}) is equivalent to 
$\displaystyle{(2\bar \gamma ^2(1-h^2)-(2\bar \beta c-b))h=0}$. So either $h=0$ and $(c+b)^2=\bar \beta ^2+\bar \gamma ^2+b$, or  $2\bar \beta c=2\bar \gamma ^2(1-h^2)+b$ and $c^2+\bar \gamma ^2(1-h^2)^2=0$. In the first case, 
 $c$ is arbitrary, depending on $b$. In the second case $c=0=\bar \gamma $, and $h$ with $\sigma (1-h^2)>0$ is arbitrary. $\Box$

\begin{example} Let $\bar \beta $ and $\bar \gamma $ be constant. If $\bar \beta =0=\bar \gamma$ then all constant curves $x$ satisfy (\ref{x1eq}), (\ref{x2eq}), (\ref{x3eq}). If $\bar \beta \not= 0$ and $\bar \gamma =0$ then the only constant solutions are $x(t)=\pm (0,0,1)$ and $x(t)=(\cos \psi _0,\sin \psi _0,0)$ where $\psi _0\in \R $. 
If neither $\bar \beta = 0$ nor $\bar \gamma = 0$,   
the only constant solutions are $x(t)=\pm (0,0,1)$. $\Box$
\end{example}

\spp
Substituting for $\Vert x^{(1)}\Vert ^2$ and $w$ in (\ref{x3eq}), we find 
$$x_3^{(2)}+bx_3~=~-2\bar \gamma ^2(1-x_3^2)x_3+\bar \gamma \bar \gamma '(1-x_3^2)^2-\bar \beta '(1-x_3^2)c+2\bar \beta x_3c$$
and, integrating with respect to $x_3$, for some $d\in \R $,
\begin{equation}\label{x3eqss}(x_3^{(1)})^2+bx_3^2~=~\bar \gamma ^2(1-x_3^2)^2-2\bar \beta c (1-x_3^2)+d.\end{equation}
More can be said for particular choices of $\bar \beta $ and $\bar \gamma$. 
\subsection{$\bar \gamma =0$ and $\bar \beta$ Constant}\label{subsecgam=0}
Let $x:\R \rightarrow M_\sigma $ be a solution of (\ref{x1eq}), (\ref{x2eq}), (\ref{x3eq}), where $\bar \gamma =0$ and $\bar \beta $ is constant. By Proposition \ref{constintegrand}, 
the integrand for $J(x,A)$ is constant. If $x_3$ is constant then $x$ is given by Proposition \ref{corconstx3}. 
\begin{theorem}\label{horthm} Let $x_3$ be nonconstant.  
For some $\lambda $ with $\sigma (1-\lambda ^2)\geq 0$ and some $v_0\in \R$, according as $\sigma =\pm 1$,  
\begin{eqnarray*}
x(t)&=&\sqrt{1-\lambda ^2\sin ^2(t\epsilon +v_0)}~(\cos \psi (t),\sin \psi (t),0)
+\lambda (0,0,\sin (t\epsilon +v_0))\quad \hbox{or}\\
x(t)&=&\sqrt{\lambda ^2\cosh ^2(t\epsilon +v_0)-1}~(\cos \psi (t),\sin \psi (t),0)+\lambda (0,0,\cosh (t\epsilon +v_0))
\end{eqnarray*}
where, for some $\sigma _1=\pm 1$,  either $\lambda ^2=1$ and $\psi (t)=t\bar \beta +\psi _0$, or $\sigma (1-\lambda ^2)>0 $ and  
\begin{eqnarray*}
\psi (t)&=&\psi _0+t\bar \beta \pm {\rm arctan} \left(  \sqrt{1-\lambda ^2}\tan (t\epsilon +v_0) \right) +2m\pi \quad \hbox{where}~~m\in \Z , 
\quad \hbox{or}\\
\psi (t)&=&\psi _0+t\bar \beta \pm {\rm arctan} \left( \frac{\tanh (t\epsilon +v_0)}{\sqrt{\lambda ^2-1}} \right) ,
\end{eqnarray*}
according as $\sigma =\pm 1$. Conversely, any such $x$ satisfies (\ref{x1eq}), (\ref{x2eq}), (\ref{x3eq}). 
\end{theorem}

\spp
{\bf Proof:}  By (\ref{x3eqss}), ~$\displaystyle{(x_3^{(1)})^2=(2\bar \beta c-b)x_3^2-2\bar \beta c+d}$. If $b=2\bar \beta c$, then $(x_3^{(1)})^2=d-2\bar \beta c$, and $d=2\bar \beta c=b$ since otherwise $x_3$ maps onto the whole of $\R $. So $x_3$ is constant when $b=2\bar \beta c$. 

\spp
For $b\not= 2\bar \beta c$, we claim that, for some $\lambda $ with $\sigma (1-\lambda ^2)\geq 0$, $\displaystyle{(x_3^{(1)})^2=(b-2\bar \beta c)(\lambda ^2-x_3^2)}$ . This holds trivially if $(2\bar \beta c-b)x_3(t)^2-2\bar \beta c+d= 0$ for some $t\in \R$. If $(2\bar \beta c-b)x_3(t)^2-2\bar \beta c+d\not= 0$ for all $t$, then  
$x_3$ is strictly monotonic, and $x_3^2$ is bounded above or below according as $\sigma =\pm 1$. So $\lim_{t\rightarrow \pm \infty}x_3(t)$ exist, $\lim_{t\rightarrow \pm \infty}x_3^{(1)}(t)=0$, and  
$$0~=~\lim_{t\rightarrow \pm \infty }x_3^{(1)}(t)^2~=~(2\bar \beta c-b)\lambda _\pm ^2-2\bar \beta c+d\quad 
\hbox{where ~} \lambda _\pm ~:=~\lim_{t\rightarrow \pm \infty}x_3(t)$$
which proves the claim. 
Set $\displaystyle{\epsilon :=\sqrt{\vert b-2\bar \beta c\vert }>0}$. Then 
$\displaystyle{(x_3^{(1)})^2=(b-2c\bar \beta )(\lambda ^2-x_3^2)}$.
\begin{itemize}
\item   
Suppose $b>2\bar \beta c$. Then $x_3(t)=\pm \lambda \sin ( t\epsilon  +v_0)$ where $v_0\in \R $. If $\sigma =-1$ then $x_3(t)^2\geq 1$ for all $t$. So $\sigma =1$ and $\lambda ^2\leq 1$. 
Then 
$$\displaystyle{x(t)=\sqrt{1-\lambda ^2\sin ^2(t\epsilon +v_0)}~(\cos \psi (t),\sin \psi (t),0)\pm \lambda (0,0,\sin (t\epsilon +v_0))}.$$
Differentiating, and substituting for $\psi ^{(1)}$ from equation (\ref{psieq}), 
\begin{eqnarray*}
\Vert x^{(1)}\Vert ^2&=&\frac{\epsilon ^2\lambda ^4\sin ^2(t\epsilon +v_0)\cos  ^2(t\epsilon +v_0)}{1-\lambda ^2\sin ^2(t\epsilon +v_0)}+\epsilon ^2\lambda ^2\cos ^2(t\epsilon +v_0)+\\
&~&(1-\lambda ^2)\sin ^2(t\epsilon +v_0)(\bar \beta +\frac{c}{1-\lambda ^2\sin ^2(t\epsilon +v_0)})^2.
\end{eqnarray*}
Equating this with the right hand side of (\ref{conssph}), $c^2=\epsilon ^2(1-\lambda ^2)$. Then the expressions for $\psi $ follow on integration of  (\ref{psieq}).  Substituting for $x$ in 
(\ref{x1eq}), (\ref{x2eq}), (\ref{x3eq}), these equations are indeed satisfied.
\item Suppose If  $b<2\bar \beta c$.  Then, for some $v_0\in \R $,  
$\displaystyle{x_3(t)=\pm \lambda \cosh (t\epsilon +v_0)}$. If $\sigma =1$ then $x_3(t)^2\leq 1$ for all $t$. So $\sigma =-1$, 
$\lambda ^2\geq 1$, and 
$$\displaystyle{x(t)=\sqrt{\lambda ^2\cosh ^2(t\epsilon +v_0)-1}~(\cos \psi (t),\sin \psi (t),0)\pm \lambda (0,0,\cosh (t\epsilon +v_0))}.$$
Proceeding as before, $c^2=\epsilon ^2(\lambda ^2-1)$,  and (\ref{x1eq}), (\ref{x2eq}), 
(\ref{x3eq}) are verified on integration of (\ref{psieq}). 
\end{itemize} $\Box$

\spp
\begin{example}\label{hor2ex} Figure \ref{horconfig} shows (red) $x:[0,1]\rightarrow S^2$ 
satisfying (\ref{x1eq}), (\ref{x2eq}), (\ref{x3eq}) with 
$\sigma =1$, $\bar \beta =-1$, $\bar \gamma =0$, $\lambda =1$, $\epsilon =7$, $v_0=\psi _0=0$.  
Figure \ref{horconfig} also shows (black) the orbits of $A$, and (blue) a numerically calculated conditional minimum 
$x_o:[0,1]\rightarrow S^2$, with 
$$x_o(0)~=~x(0)~=~(1,0,0)\quad \hbox{and}\quad x_o(1)~=~x(1)~=~(\cos 1\cos 7,-\cos 7\sin 1,\sin 7).$$
Then  $J(x_o,A)\approx 0.519$ and $J(x,A)=49$. Repeating with $\epsilon =\sqrt{2}$ gives $x_o\approx x$,   
%
$J(x_o,A)\approx 2.011$ and $J(x)=2$.  
\begin{figure}[htbp] 
   \centering
   \includegraphics[width=2in]{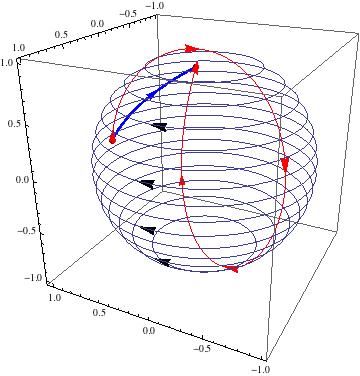} 
   \caption{A Non-Minimal Conditional Extremum (red), with the Minimum (blue), in Example \ref{hor2ex}}
   \label{horconfig}
\end{figure}

\spp
Figure \ref{hor1fig} shows $x:[0,10]\rightarrow S^2$ with $\sigma =1$, $\bar \beta =-1$, $\bar \gamma =0$, $\lambda =1$. 
\begin{figure}[htbp] 
   \centering
   \includegraphics[width=2in]{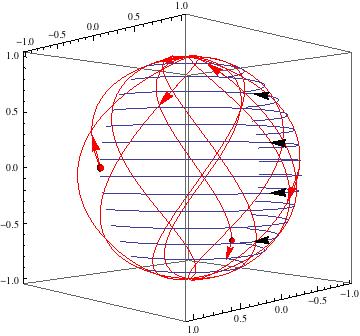} 
   \caption{$\sigma =1$, $\bar \beta =-1$, $\lambda =1$, $\epsilon ^2=1.1$, $v_0=\psi _0=0$ in Example \ref{hor2ex}}
   \label{hor1fig}
\end{figure}

\spp
Figure \ref{hor2fig} shows $x:[0,10]\rightarrow S^2$ with $\sigma =1$, $\bar \beta =-1$, $\bar \gamma =0$, $\lambda =0.5$, $\epsilon =1.1$, $v_0=\psi _0=0$. 
\begin{figure}[htbp] 
   \centering
   \includegraphics[width=2in]{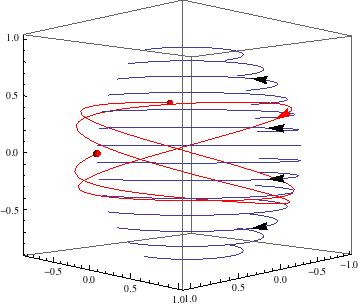} 
   \caption{$\sigma =1$, $\bar \beta =-1$, $\lambda =0.5$, $\epsilon =1.1$, $v_0=\psi _0=0$ in Example \ref{hor2ex}}
   \label{hor2fig}
\end{figure}
$\Box$
\end{example}

\begin{example}\label{poin1ex} The Poincar\'e disc $D$ is the open unit disc in $\R ^2$ with Riemannian metric  
$$\langle v,w\rangle ~:=~4\frac{(v_1w_1+v_2w_2)}{(1-z_1^2+z_2^2)^2}$$ 
where $v,w$ are tangent to $D$ at $z$. Geodesics of $D$ are either circular arcs orthogonal to the unit circle bounding $D$, or diameters. An isometry to $D$ from the path component of  $(0,0,1)$ in $H^2$ is given by 
$$y~\mapsto ~\phi (y)~:=~\frac{(y_1,y_2)}{1+y_3}.$$ 
Rather than directly plot critical curves in the hyperboloid $H^2$, a clearer picture is obtained by mapping to $D$. Figure \ref{poin1fig} shows $\phi \circ x:[0,10]\rightarrow D$ (red) with $\sigma =-1$, $\bar \beta =-1$, $\lambda =1.01$, $\epsilon=4$, $v_0=0$ and $\psi _0=-\pi /3$. Integral curves of the field corresponding to $A$ are shown as circles (black). Where the conformal factor in the metric is large, towards the end, the conditional minimum $x$ tends to agree with the orbits of the prior field. Disagreements with the prior field are greater towards the beginning. 
%
\begin{figure}[htbp] 
   \centering
   \includegraphics[width=2in]{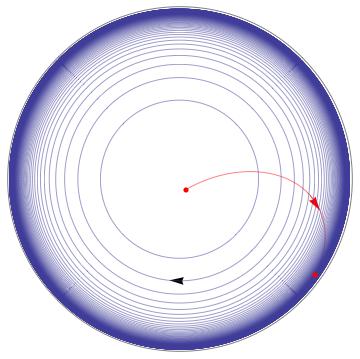} 
   \caption{$\sigma =-1$, $\bar \beta =-1$, $\lambda =1.01$, $\epsilon=4$, $v_0=0$, $\psi _0=-\pi /3$ in Example \ref{poin1ex}}
   \label{poin1fig}
\end{figure}
$\Box$
\end{example}
\subsection{$\bar \beta $, $\bar \gamma $ Constant and $\bar \gamma \not=  0$}\label{gamnon0sec} 
Let $x:\R \rightarrow M_\sigma$ be a solution of (\ref{x1eq}), (\ref{x2eq}), (\ref{x3eq}), where $\bar \beta $ and $\bar \gamma $ are constant with $\bar \gamma \not= 0$. By Theorem \ref{thm2}, $A$ is reflexive if and only if $\bar \beta =0$. 
If $x_3$ is constant then $x$ is given by part 1 of Proposition \ref{corconstx3}. 
\begin{theorem}\label{wthm} If $x_3$ is nonconstant then, for some $a\in \C $ 
$$x_3(t)~=~\pm \sqrt{\wp (\bar \gamma t+a;g_2,g_3)-\bar \delta }$$
where $\bar \delta :=(2\bar \beta c-b-2\bar \gamma ^2)/(3\bar \gamma ^2)$, $\wp $ is the Weierstrass elliptic function and, for some $\bar d\in \R $, 
\begin{eqnarray*}
g_2&=&~~12\bar \delta ^2+\bar d \\
g_3&=&-8\bar \delta ^3-\bar \delta \bar d .
\end{eqnarray*}
\end{theorem}

\spp
{\bf Proof:} By (\ref{x3eqss}),  
$\displaystyle{(x_3^{(1)})^2=\bar \gamma ^2(x_3^4+3\bar \delta x_3^2-\bar d/4)}$
where $\bar d:=-4(\bar \gamma ^2-2\bar \beta c+d)/\bar \gamma ^2$. Setting $\displaystyle{y(s):=x_3(s/\bar \gamma )^2+\bar \delta }$, we find  
$\displaystyle{\left( \frac{dy}{ds}\right) ^2=4y^3-g_2y-g_3}$. $\Box$

\spp
So (\ref{psieq}) gives $x$ up to quadrature in terms of a known periodic function $x_3:\R \rightarrow \R$. 
\begin{example}\label{lamex} Take $\sigma =1$, $\bar \beta =0$ and $\bar \gamma =1$. Since $A$ is conservative it is reflexive, by 
Theorem \ref{thm2}, and the reverse $\bar x$ of a conditional extremum is also a conditional extremum for the reversed data. 
%
%
Figure \ref{fig0} shows (red) the solution $x:[0,14]\rightarrow S^2$ of 
(\ref{x1eq}), (\ref{x2eq}), (\ref{x3eq}) with $x(0)=x_0=(\sqrt{3}/2,0,1/2)$ and $x^{(1)}(0)=(0,1,0)$. The integral curves of 
$A$, namely the gradient of the height function, are shown in black.  
\begin{figure}[htbp] 
   \centering
   \includegraphics[width=2in]{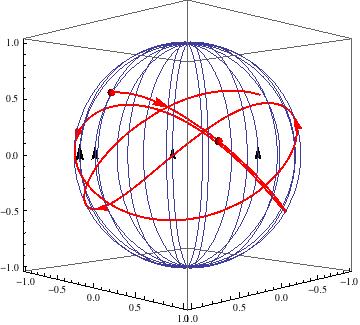} 
   \caption{A Conditional Extremum for the Conservative Prior Field in Example \ref{lamex}}
   \label{fig0}
\end{figure}
%
We find  $a\approx 1.14811+1.74899 \imath$, 
$b=0.25$, $c=-0.5$, $g_2=4.75$, $g_3=1.875$. $\Box$  
\end{example}

\begin{example}\label{poin2ex} Take $\sigma =-1$, $\bar \beta =-1$ and 
$\bar \gamma =2$. Figure \ref{poin2fig} shows (red) $\phi \circ x:[0,1]\rightarrow D$ 
where $x:[0,1]\rightarrow H^2$ is the conditional extremum with 
$x(0)=x_0=(0.1,0.1,\sqrt{1.02})$, $x(1)=x_1=(8.99009,-7.34992,11.6552)$ and $x^{(1)}(0)=(0.9\sqrt{1.002},0,0.9)$. Integral curves of the prior field $A$ spiral inwards (black). Numerical calculations suggest 
$x$ is a conditional minimum. 
\begin{figure}[htbp] 
   \centering
   \includegraphics[width=2in]{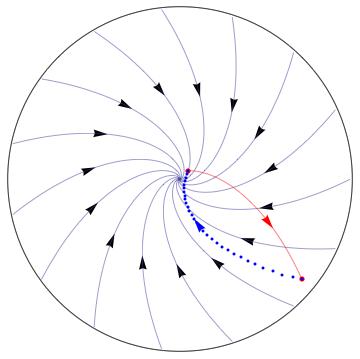} 
   \caption{Conditional Minima for the NonConservative Prior Field in Example \ref{poin2ex}}
   \label{poin2fig}
\end{figure}
We find $g_2=1.17393$, $g_3 =-0.220814$ and $a\approx -2.08599$. The prior field $A$ is effectively contradicted by the supposed observations $x_0$ and $x_1$. The conditional minimum $x$ makes the best of difficult circumstances. 

\spp
On the other hand, if the order of $x_0$ and $x_1$ is interchanged, these might well lie on an integral curve of a vector field near $A$. Figure \ref{poin2fig} shows the corresponding conditional minimum (blue, dotted). Disagreements with $A$ tend to be concentrated near the centre of the disc, where the conformal factor in the Riemannian metric is relatively small. By Theorem \ref{thm2}, $A$ is not reflexive, because $\bar \beta \not= 0$. 
$\Box$
\end{example}

\spp
{\bf Acknowledgement:} This paper is dedicated, with great respect, to the memory of Professor Jerrold E. Marsden.

\end{document}